\documentclass[12pt]{article}

\usepackage{amsfonts,amssymb,amsmath,amsthm,epsfig}

\usepackage{graphics,graphicx}                 
\usepackage{xypic}

\def\boxit#1#2{\hbox{\vrule
    \vtop{%
      \vbox{\hrule\kern#1%
        \hbox{\kern#1#2\kern#1}}%
      \kern#1\hrule}%
    \vrule}}

\newtheorem{theorem}{Theorem}

\theoremstyle{definition}
{

\newtheorem{problem}{Problem}

}

\title{Crucial and bicrucial permutations with respect to arithmetic monotone patterns}

\author{
Sergey Avgustinovich\\
\small Sobolev Institute of Mathematics\\[-0.8ex]
\small 4 Acad. Koptyug Ave\\[-0.8ex]
\small 630090 Novosibirsk, Russia\\[-0.8ex]
\small \texttt{avgust@math.nsc.ru}
\and
Sergey Kitaev \\
\small Department of Computer and Information Sciences\\[-0.8ex]
\small University of Strathclyde\\[-0.8ex]
\small Glasgow G1 1XH, United Kingdom\\[-0.8ex]
\small \texttt{sergey.kitaev@cis.strath.ac.uk}
\and
Alexandr Valyuzhenich \\
\small Novosibirsk State University\\[-0.8ex]
\small 2 Pirogova Street\\[-0.8ex]
\small 630090 Novosibirsk, Russia\\[-0.8ex]
\small \texttt{graphkiper@mail.ru}
}

\begin{document}
\maketitle

\begin{abstract}
\noindent \

A pattern $\tau$ is a permutation, and an arithmetic occurrence of $\tau$ in (another) permutation $\pi=\pi_1\pi_2\ldots\pi_n$ is a subsequence $\pi_{i_1}\pi_{i_2}\ldots\pi_{i_m}$ of $\pi$ that is order isomorphic to $\tau$ where the numbers $i_1<i_2<\cdots<i_m$ form an arithmetic progression. A permutation is $(k,\ell)$-crucial if it avoids arithmetically the patterns $12\ldots k$ and $\ell(\ell-1)\ldots 1$ but its extension to the right by any element does not avoid arithmetically these patterns. A $(k,\ell)$-crucial permutation that cannot be extended to the left without creating an arithmetic occurrence of $12\ldots k$ or $\ell(\ell-1)\ldots 1$ is called $(k,\ell)$-bicrucial. 

In this paper we prove that arbitrary long $(k,\ell)$-crucial and $(k,\ell)$-bicrucial permutations exist for any $k,\ell\geq 3$. Moreover, we show that the minimal length of a $(k,\ell)$-crucial permutation is $\max(k,\ell)(\min(k,\ell)-1)$, while the minimal length of a $(k,\ell)$-bicrucial permutation is at most $2\max(k,\ell)(\min(k,\ell)-1)$, again for $k,\ell\geq3$.\\


\end{abstract}

\section{Introduction}

The notion of a {\em crucial word} with respect to a given set of prohibitions is introduced in \cite{EK}. A word is crucial if it avoids a given set of prohibitions but any extension of the word to the right by a letter of the alphabet in question does not avoid the prohibitions. Crucial words were studied with respect to several sets of prohibitions in the literature \cite{AGHK,EK,GHK}. A crucial word is called {\em bicrucial} if it is also not possible to extend it to the left without creating a prohibition. Bicrucial words are called {\em maximal words} in the literature and they were studied with respect to so-called {\em squares}  and {\em abelian squares} \cite{C,K}. An example of a (bi)crucial word over the alphabet $\{1,2,\ldots,n\}$ with respect to squares/abelian squares is the {\em Zimin word} $X_n$ defined recursively as $X_1=1$ and $X_n=X_{n-1}nX_{n-1}$.

An extension of the notion of a crucial word to a crucial permutation \cite{AKPV} was closely related to the {\em theory of patterns in permutations and words} (see \cite{Kit} for a comprehensive introduction to the area). A pattern is a permutation or word $\tau$ involving each letter from an alphabet $\{1,2,\ldots,k\}$ for some $k$, and an occurrence of $\tau$ in another permutation or word $\pi$ is a subsequence of $\pi$ that is order isomorphic to $\tau$. For a permutation $\pi$, its {\em reduced form} is the permutation obtained from $\pi$ by substituting the $i$th largest element by $i$. For example, the reduced form of 2754 is 1432. Thus, an occurrence of $\tau$ in $\pi$ is a subsequence of $\pi$ whose reduced form is $\tau$. When occurrences of $\tau$ are required to be formed by consecutive elements, $\tau$ is called a {\em consecutive pattern}. If there are two occurrences of a consecutive pattern of length at least 2 in a permutation $\pi$ that are next to each other, we say that $\pi$ contains a square. (Bi)crucial permutations with respect to squares were studied in \cite{AKPV}, while the present paper deals with {\em arithmetic occurrences} of monotone patterns in permutations and (bi)crucial permutations associated with them (see below for precise definitions). 

The celebrated {\em Van der Waerden's theorem} states that for any given positive integers $r$ and $k$, there is some number $N$ such that if the integers $\{1,2,\ldots,N\}$ are colored, each with one of $r$ different colors, then there will be at least $k$ integers forming an arithmetic progression all of the same color. An intriguing fact is that it is impossible to extend the van der Waerden's theorem to the case of infinite permutations. Partially, it is explained by existence of so-called anti-monotone permutations in which no three elements are increasing or decreasing in positions forming an arithmetic progression; see \cite{DEGS}. However, it is interesting to study the structure of arithmetic progressions in  infinite permutations. In particular, Makarov \cite{Mak,Mak1} studied {\em arithmetic complexity}, defined similarly to the notion of arithmetic complexity of words \cite{AFF,Frid1}, for a number of well-known infinite permutations. For some relevant work, see \cite{AFKS,FFF,Frid}.

Let $\pi=\pi_1\pi_2\ldots\pi_n$ be a permutation. Then for a fixed $d\geq 1$, $\pi_i\pi_{i+d}\ldots\pi_{i+(k-1)d}$ is an {\em arithmetic subsequence} of length $k$ with difference $d$, assuming $i\geq 1$ and $i+(k-1)d\leq n$. If an occurrence of a pattern form an arithmetic subsequence, we refer to the occurrence as an {\em arithmetic occurrence} of the pattern. A permutation $\pi$ {\em avoids arithmetically} a pattern $\tau$ if $\pi$ contains no arithmetic occurrences of $\tau$. A permutation is {\em $(k,\ell)$-anti-monotone} if it avoids arithmetically the patterns $12\ldots k$ and $\ell(\ell-1)\ldots 1$. See \cite{Sharma} for relevant work, where the operation of taking (the usual group-theoretic) inverse can be applied to all permutations to enter our domain.

A useful fact appearing in \cite{DEGS}, to be used throughout the paper, is that (3,3)-anti-monotone permutations exist of arbitrary length. For example, a way to construct a (3,3)-anti-monotone permutation of length $2n$ is to take two (3,3)-anti-monotone permutations of length $n$, say $\pi_1\pi_2\ldots\pi_n$ and $\tau_1\tau_2\ldots\tau_n$ and to consider the {\em shuffle down} of these permutations $(\pi_1+n)\tau_1(\pi_2+n)\tau_2\ldots(\pi_n+n)\tau_n$ to get the sought permutation; taking two (3,3)-anti-monotone permutations $\pi_1\pi_2\ldots\pi_n$ and $\tau_1\tau_2\ldots\tau_{n-1}$ would give a (3,3)-anti-monotone permutation of length $2n-1$, say,  $(\pi_1+n-1)\tau_1(\pi_2+n-1)\tau_2\ldots\tau_{n-1}(\pi_n+n-1)$.

An {\em extension} of a permutation $\pi$ of length $n$ to the right (resp., left) is a permutation $\pi'x$ (resp., $x\pi'$) of length $n+1$ such that $x\in\{1,2,\ldots,n+1\}$ and $\pi'$ is obtained from $\pi$ by adding 1 to each element that is large than or equal to $x$.  For example, the set of all extensions to the right of the permutation $231$ is $\{3421,3412,2413,2314\}$, while the set of all extensions to the left of $12$ is $\{123,213,312\}$.

A permutation $\pi$ is {\em $(k,\ell)$-anti-monotone crucial} (resp. {\em $(k,\ell)$-anti-monotone bicrucial}) if $\pi$ is anti-monotone but any extension of $\pi$ to the right (resp., and to the left) is not  $(k,\ell)$-anti-monotone. $(k,\ell)$-anti-monotone crucial permutations and $(k,\ell)$-anti-monotone bicrucial permutations will be called by us {\em $(k,\ell)$-crucial} and {\em $(k,\ell)$-bicrucial} permutations, respectively, for brevity.  For example, the permutation 216453 is $(3,3)$-crucial, while the permutation 73418562 is $(3,3)$-bicrucial. 

If $\pi=\pi_1\pi_2\ldots\pi_n$ is a permutation, then its {\em reverse}, $r(\pi)$, is the permutation $\pi_n\pi_{n-1}\ldots \pi_1$, while the {\em complement} of $\pi$, $c(\pi)$, is obtained from $\pi$ by replacing an element $i$ by $n+1-i$, for $1\leq i\leq n$. For example, if $\pi=24135$ then $r(\pi)=53142$ and $c(\pi)=42531$.

Clearly, if $\pi$ is a $(k,\ell)$-(bi)crucial permutation then $c(\pi)$ is an $(\ell,k)$-(bi)crucial permutation. Thus in our studies, without loss of generality, we can assume that $k\leq \ell$. Another straightforward fact is that if $\pi$ is a $(k,\ell)$-bicrucial permutation then $r(\pi)$ is an $(\ell,k)$-bicrucial permutation.

We are interested in the following questions in which we assume that $k,\ell\geq 3$: Do there exist $(k,\ell)$-(bi)crucial permutations? If so, what is the minimum length of such permutations? Do arbitrary long $(k,\ell)$-(bi)crucial permutations exist? In this paper we will prove that arbitrary long $(k,\ell)$-crucial and  $(k,\ell)$-bicrucial permutations exist for any $k,\ell\geq 3$. Moreover, we will show that the minimum length of a $(k,\ell)$-crucial permutation, denoted $m(k,\ell)$ is $\max (k,\ell)(\min (k,\ell)-1)$, while the minimum length of a $(k,\ell)$-bicrucial permutation is bounded from above by $2m(k,\ell)$. $(k,\ell)$-(bi)crucial permutations of smallest length are called {\em minimal} $(k,\ell)$-bicrucial permutations.

The paper is organized as follows. In Section \ref{kl-crucial-lower} we give a lower bound for $m(k,\ell)$, while an upper bound for $m(k,\ell)$ is discussed in Section \ref{kl-crucial-upper}. In Section \ref{arbitraryLong} we discuss existence of arbitrary long $(k,\ell)$-crucial permutations and in Section \ref{kl-bicrucial}  we discuss $(k,\ell)$-bicrucial permutations. Finally, in Section \ref{final} we state some open problems. 

\section{A lower bound for the minimal length of a $(k,\ell)$-crucial permutation}\label{kl-crucial-lower}

In this section we will prove that for $k,\ell\geq 3$, $m(k,\ell)$, the minimum length of  a $(k,\ell)$-crucial permutation is $\geq\max(k,\ell)(\min(k,\ell)-1)$.

Let $\pi=\pi_1\pi_2\ldots\pi_n$ be a $(k,\ell)$-crucial permutation of minimal length, that is, $n=m(k,\ell)$. If a subsequence $\pi_{i_1}\pi_{i_2}\ldots\pi_{i_m}$ is an arithmetic occurrence of the pattern $12\ldots m$ (resp., $m(m-1)\ldots 1$) we say that the subsequence is an {\em AP($m$)-up permutation} (resp., {\em AP($m$)-down permutation}). Further we let $M$ (resp., $N$) be the set of all AP($k-1$)-up permutations (resp., AP($\ell$)-down permutations) of the form   $\pi_{n+1-(k-1)a}\pi_{n+1-(k-2)a}\ldots\pi_{n+1-a}$ (resp., $\pi_{n+1-(\ell-1)b}\pi_{n+1-(\ell-2)b}\ldots\pi_{n+1-b}$) for some $a,b\geq 1$. We have $M\neq\emptyset$ because extending $\pi$ to the right by $n+1$ we must get a prohibited monotone subsequence, but $n+1$ cannot be the last element in a decreasing subsequence of length $\ell$. Similarly,  $N\neq\emptyset$ because extending $\pi$ to the right by the smallest element $1$ we must get a prohibited monotone subsequence, but $1$ cannot be the last element in an increasing subsequence of length $k$. 

We let $a^*$ (resp. $b^*$) be the minimum (resp., maximum) possible rightmost element in a permutation in $M$ (resp., in $N$). Next we will prove that $a^*<b^*$. 

Suppose $a^*>b^*$ and consider extending $\pi$ to the right by an element $x\in\{b^*+1,\ldots,a^*\}$. We claim that the obtained permutation, say $\pi'$, is $(k,\ell)$-anti-monotone which contradicts to the fact that $\pi$ is a $(k,\ell)$-crucial permutation. Indeed, if $\pi'$ contains an arithmetic occurrence of the pattern $12\ldots k$ ending with $x$, then next to last element, say $y$, in this occurrence must be less than $a^*$ (the element $y$ of $\pi$ was unchanged in $\pi'$), a contradiction with the definition of $a^*$. On the other hand, if $\pi'$ contains an arithmetic occurrence of the pattern $\ell(\ell-1)\ldots 1$ ending with $x$, then next to last element, say $z$, in this occurrence must be larger than $b^*+1$ in $\pi'$ and thus $z$ is larger than $b^*$ in $\pi$,  a contradiction with the definition of $b^*$. Thus, we must have $a^*<b^*$. 

Suppose now that an AP($k-1$)-up permutation  $\pi_{n+1-(k-1)a}\pi_{n+1-(k-2)a}\ldots\pi_{n+1-a}$ is such that $\pi_{n+1-a}=a^*$ and an AP($\ell-1$)-down permutation  $\pi_{n+1-(\ell-1)b}\pi_{n+1-(\ell-2)b}\ldots\pi_{n+1-b}$ is such that $\pi_{n+1-b}=b^*$. Since $a^*<b^*$, these permutations have no elements in common and therefore the positions of all the elements are different.  Thus, we have two sets $P=\{a,2a,\ldots,(k-1)a\}$ and $Q=\{b,2b,\ldots,(\ell-1)b\}$ such that $a,b\geq 1$ and $P\cap Q=\emptyset$.

Suppose that $(k-1)a<(\ell-1)b$ then $ba\not\in P$ (since otherwise $ba=ab\in Q$, a contradiction with $P\cap Q=\emptyset$). Thus $b\geq k$ and $(\ell-1)b\geq (\ell-1)k$, from where the length $n$ of $\pi$ is no less than $(\ell-1)k$. 

Suppose now that $(\ell-1)b<(k-1)a$ then $ab\not\in Q$ (since otherwise $ab=ba\in P$, a contradiction with  $P\cap Q=\emptyset$). Thus $a\geq\ell$ and $(k-1)a\geq (k-1)\ell$, from where the length $n$ of $\pi$ is no less than $(k-1)\ell$. 

Summarizing our considerations, the length $n$ of $\pi$ cannot be less than 
$$\min((\ell-1)k,(k-1)\ell)=\max(k,\ell)(\min(k,\ell)-1).$$   

\section{An upper bound for the length of a minimal $(k,\ell)$-crucial permutation}\label{kl-crucial-upper}

In this section, we consider three cases depending on whether $\min(k,\ell)\geq 5$, $\min(k,\ell)=4$ or $\min(k,\ell)= 3$ (see the respective subsections) to prove that the minimum length of a $(k,\ell)$-crucial permutation is at most $\max(k,\ell)(\min(k,\ell)-1)$. Then summarizing considerations in these subsection and Section \ref{kl-crucial-lower} we will get a proof of the following theorem.

\begin{theorem}\label{MainTh} For $k,\ell\geq 3$, the minimum length of a $(k,\ell)$-crucial permutation, $m(k,\ell)$, is equal to $\max(k,\ell)(\min(k,\ell)-1)$.\end{theorem}

\subsection{The case of $\min(k,\ell)\geq 5$.}

\begin{figure}[ht]
\begin{center}
\includegraphics[scale=0.5]{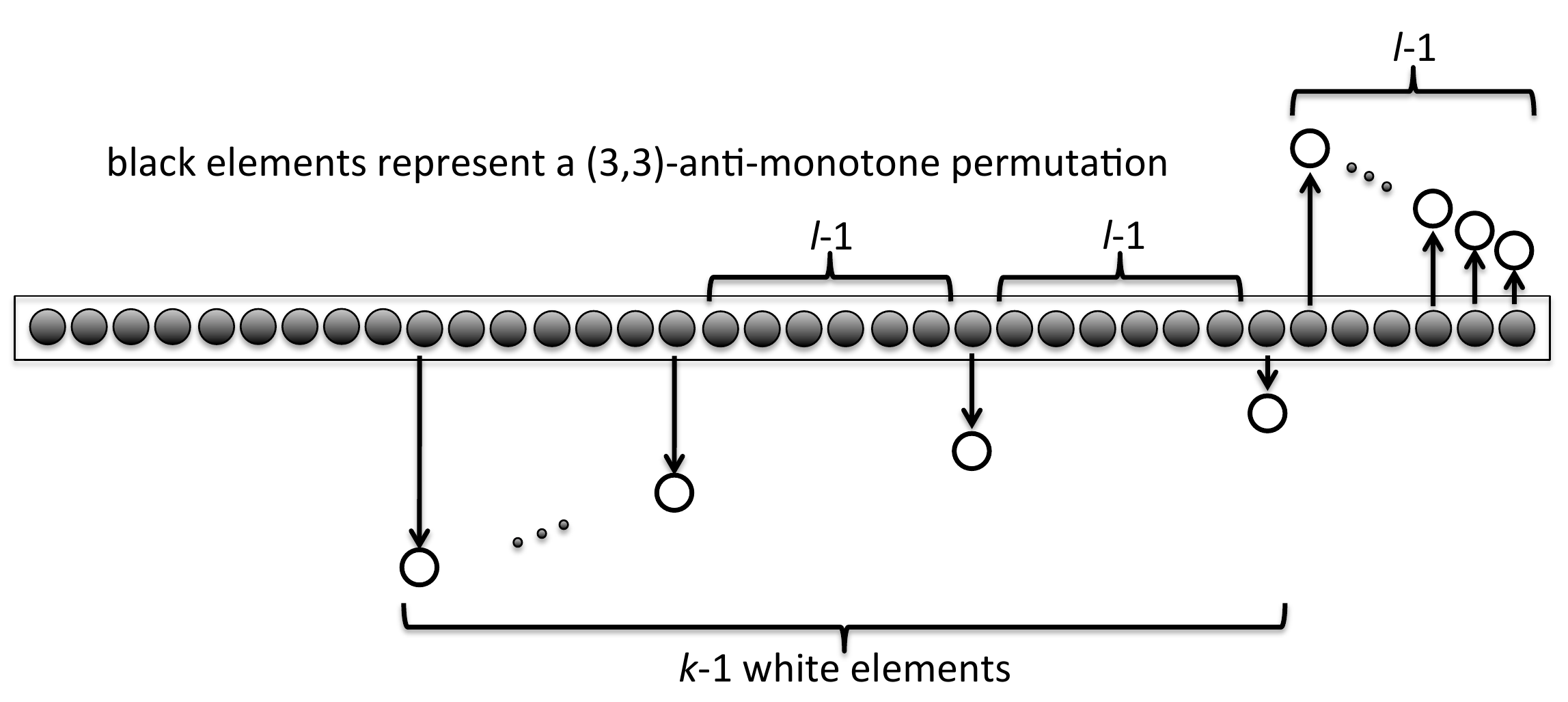}
\end{center}
\caption{A construction of $(k,\ell)$-crucial permutations, $k,\ell\geq 5$ and $k\leq \ell$.}\label{mainConstruction} 
\end{figure}

We begin with presenting a construction, shown schematically in Figure \ref{mainConstruction}, of a $(k,\ell)$-crucial permutation of any length greater than or equal to $\ell(k-1)$ for $k,\ell\geq 5$ and $k\leq\ell$. Let $\pi$ be an arbitrary $(3,3)$-anti-monotone permutation of length $n$, $n\geq \ell(k-1)$, shown schematically with black circles in Figure \ref{mainConstruction}. Rewrite the elements of $\pi$  in positions $1, 2,\ldots,k-1$ counted from right to form the permutation $n(n-1)\ldots (n-\ell+2)$ read from left to right, and rewrite the subsequence in positions $k,2k,\ldots,(\ell-1)k$ counted from right to form the permutation $12\ldots(k-1)$ read from left to right; keep the relative order of the remaining elements of $\pi$ the same, and let them form a permutation on the set $\{k,k+1,\ldots,n-\ell+1\}$. For example,  let $\pi=1(17)9(13)5(15)7(11)3(16)8(12)4(14)6(10)2$ be a $(3,3)$-anti-monotone permutation and $k=4$ and $\ell=5$. Then the permutation obtained by applying the construction in Figure \ref{mainConstruction} is 4(13)1(10)6(11)725(12)893(17)(16)(15)(14).

It is straightforward to prove that the construction is valid, namely that it presents a $(k,\ell)$-anti-monotone permutation, say $\pi'$, that is not extendable to the right. Indeed, suppose $\pi'$ contains an arithmetic occurrence of  the pattern $\ell(\ell-1)\ldots 1$. Clearly, such an occurrence cannot involve the $\ell$ rightmost elements of $\pi'$ (which are $\{n-\ell+2,n-\ell+3,\ldots,n\}$) and it can involve at most one small white element (one element from the set $\{1,2,\ldots,k\}$).  However, since only two black elements can be involved in that occurrence (the black elements form a $(3,3)$-anti-monotone permutation) and $\ell\geq 5$ (for that part, we actually need $\ell\geq 4$), no such occurrence is possible. Similarly, the longest arithmetic occurrence of an increasing pattern can involve at most two black elements, at most one large white element and at most one small white element, and thus such an occurrence is not possible in $\pi'$ since $k\geq 5$. Finally, trying to extend $\pi'$ to the right by an element larger than $k-1$, we will get an arithmetic occurrence of the pattern $12\ldots k$ (involving the $k-1$ smallest elements in $\pi'$), while extending it by an element smaller than $k$, we will get an arithmetic occurrence of the pattern $\ell(\ell-1)\ldots1$ (involving the $\ell-1$ rightmost elements in $\pi'$).  

Applying the complement operation to the construction in Figure \ref{mainConstruction} and switching $k$ and $\ell$ we cover the case of $k,\ell\geq 5$ and $k\geq \ell$. In either case, when $\min(k,\ell)\geq 5$, we have that $m(k,\ell)\leq \max(k,\ell)(\min(k,\ell)-1)$.

\subsection{The case of $\min(k,\ell)=4$.}

As the matter of fact, the construction in  Figure \ref{mainConstruction} is suitable for the case $k=\ell=4$. The only difference is that there is one restriction on the (3,3)-anti-monotone permutation $\pi$: the 6th element from the right must be smaller than the 9th element from the right in order to avoid an arithmetic occurrence with difference $d=3$ of the pattern $1234$ involving the largest element and the smallest element.  An example of a (4,4)-crucial permutation obtained in this way is $\pi=185926743(12)(11)(10)$; arbitrary long such permutations can also be constructed using an appropriate shuffling for creating (3,3)-anti-monotone permutations, but we do not provide more details here. Note that the length of $\pi$ is $12=\max(k,\ell)(\min(k,\ell)-1)$.

\begin{figure}[ht]
\begin{center}
\includegraphics[scale=0.5]{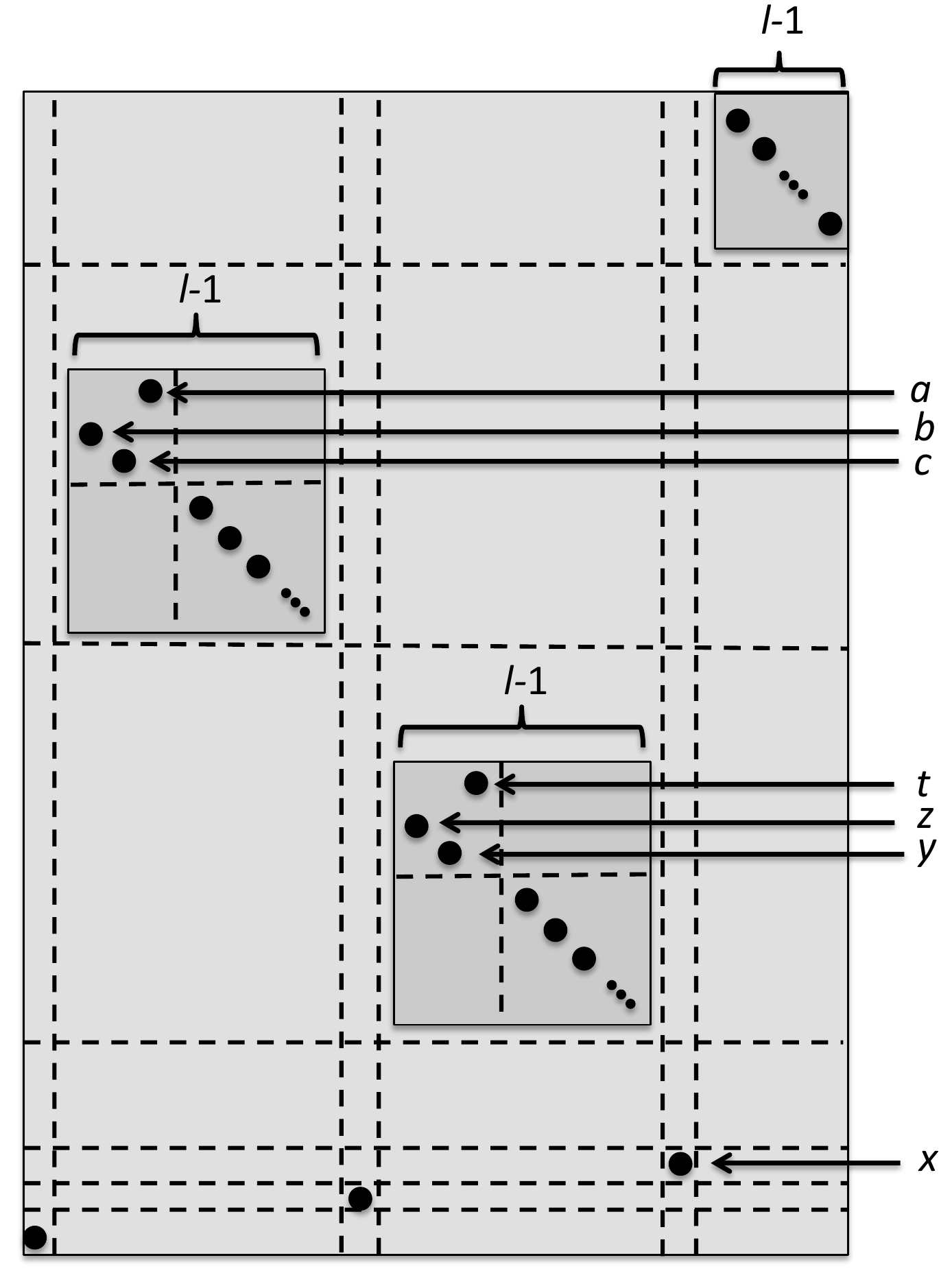}
\end{center}
\caption{A construction of a $(4,\ell)$-crucial permutation, $\ell\geq 4$.}\label{4l-construction} 
\end{figure}

In either case, we need a construction for  $(4,\ell)$-crucial permutations where $\ell\geq 5$. We provide such a construction (that actually works for $\ell=4$ as well, and gives the desired length $3\ell$) in Figure \ref{4l-construction}. Examples of permutations produced by the construction are $187925463(12)(11)(10)$ for $\ell=4$ and 
$1  (10)9(11)8 2  6574 3 (15)(14)(13)(12)$  for $\ell=5$. It is easy to see that the permutation, whose structure is presented in Figure \ref{4l-construction}, is not extendable to the right. We only need to justify that this permutation, say $\sigma$, avoids arithmetically the patterns $1234$ and $\ell(\ell-1)\ldots 1$. Indeed, if $\sigma$ contains arithmetically $1234$ then, as it is easy to see, exactly one of the elements 1, 2 and 3 must be involved, and exactly one of the $\ell-1$ largest elements must be involved, which means that the remaining middle elements in the occurrence of 1234 are either $ba$ or $ca$ or $zt$ or $yt$; it is straightforward to check that none of the four cases is possible. On the other hand, to get an arithmetic occurrence of the pattern $\ell(\ell-1)\ldots 1$ in $\sigma$, the difference of the arithmetic progression must be $d=2$ (the case $d=1$ is straightforward to see not to work, while for $d>2$ we clearly will not get enough elements since the $\ell-1$ topmost elements cannot be involved -- the rightmost element cannot be to the right of $x$). However, if $d=2$ then any arithmetic occurrence of $\ell(\ell-1)\ldots 1$ will necessarily involve either $2y$ or $zt$ which is impossible.

Applying the complement operation to Figure \ref{4l-construction}, we take care of the case of $(k,4)$-crucial permutations. Thus, in the case of $\min(k,\ell)=4$, we have that $m(k,\ell)\leq \max(k,\ell)(\min(k,\ell)-1)$.

\subsection{The case of $\min(k,\ell)=3$.}

As it was mentioned above,  the permutation 216453 is $(3,3)$-crucial and it is of length $6=\max(k,\ell)(\min(k,\ell)-1)$. So we only need to consider the case of $(3,\ell)$-crucial permutations with $\ell\geq 4$ (then applying the complement operation we will take care of the case of $(k,3)$-crucial permutations with $k\geq 4$).

\begin{figure}[ht]
\begin{center}
\includegraphics[scale=0.5]{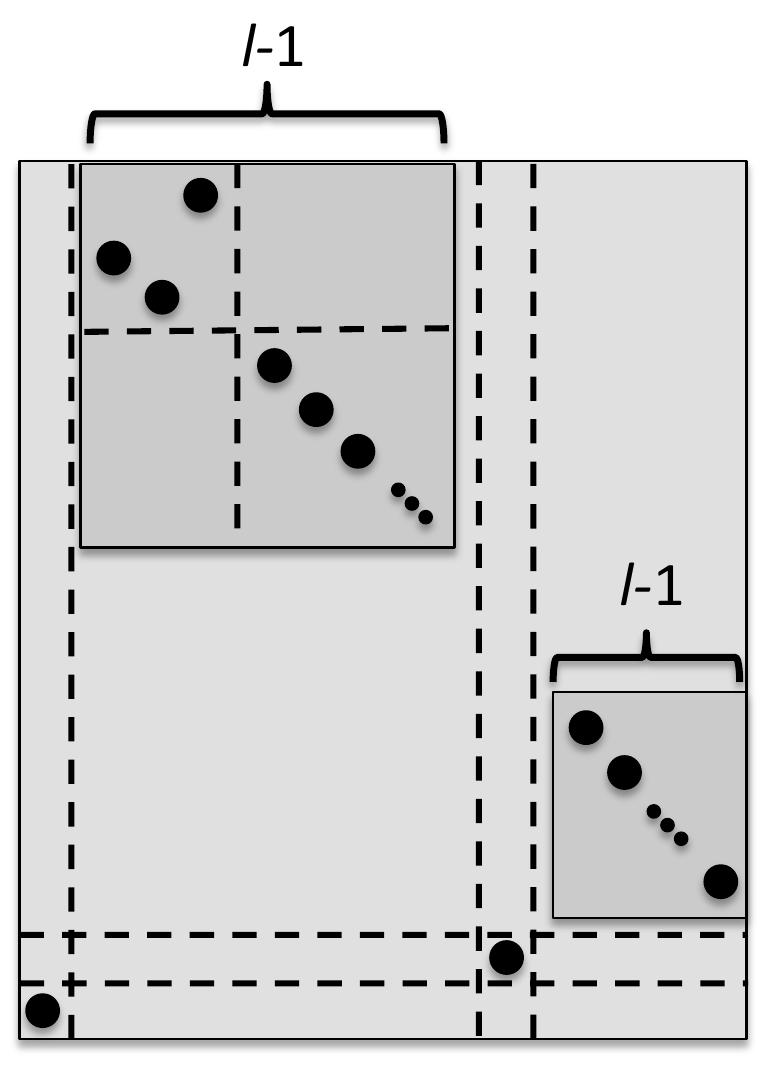}
\end{center}
\caption{A construction of a $(3,\ell)$-crucial permutation, $\ell\geq 4$.}\label{3l-construction} 
\end{figure}

It is straightforward to see that the permutation, say $\pi$, constructed in Figure \ref{3l-construction} avoids arithmetically the pattern $123$ and it has no arithmetic occurrences of the pattern $\ell(\ell-1)\ldots 1$ with differences $d=1$ and $d>2$ (the latter follows from the fact that the permutation is of length $2\ell$). Thus $d=2$ and both of the elements around 2 must be involved in the  arithmetic occurrence of $\ell(\ell-1)\ldots 1$. If $\ell=2k$ for some $k$, then the occurrence can  use at most $k$ elements to the right of 2 and at most $k-1$ elements to the left of 2 (because of the two largest elements) giving in total less than $\ell$ elements, a contradiction. On the other hand, if $\ell=2k+1$, the occurrence can use at most $k$ elements to the right of 2, and at most $k$ elements to the left of 2, thus in total less than $\ell$ elements, a contradiction.

Thus, in the case of $\min(k,\ell)=3$, we have  $m(k,\ell)\leq\max(k,\ell)(\min(k,\ell)-1)$.

\section{Arbitrary long $(k,\ell)$-crucial permutations}\label{arbitraryLong}

\begin{figure}[ht]
\begin{center}
\includegraphics[scale=0.5]{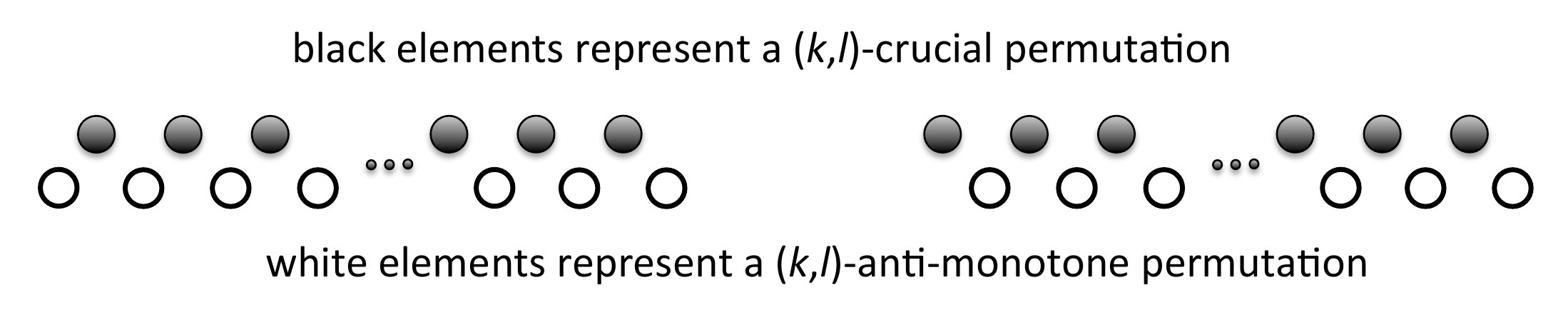}
\end{center}
\caption{Two constructions of $(k,\ell)$-crucial permutations of increasing lengths based on smaller  $(k,\ell)$-crucial permutations; $k,\ell\geq 3$.}\label{longCrucial} 
\end{figure}

From Figure \ref{mainConstruction} it is clear that arbitrary long $(k,\ell)$-crucial permutations exist of any length larger than or equal to  $\max(k,\ell)(\min(k,\ell)-1)$ for $k,\ell\geq 5$. Existence of arbitrary long $(3,3)$-crucial permutations follows, for example, from \cite{Mak1}. In either case, in Figure \ref{longCrucial} we present two constructions of obtaining longer $(k,\ell)$-crucial permutations assuming we have a construction of shorter $(k,\ell)$-crucial permutations for any $k,\ell\geq 3$, which is always the case by considerations above; however, we cannot guarantee existence of respective crucial words of any length. 

In Figure \ref{longCrucial} , the black circles represent schematically a given  $(k,\ell)$-crucial permutation, say $\pi$ of length $n$, while the white circles represent schematically an arbitrary $(k,\ell)$-anti-monotone permutation $\tau$, respectively, of length $n+1$ and $n$ (in fact, the constructions will work if we will take instead an arbitrary $(3,3)$-anti-monotone permutation). Each white element must be smaller than every black element. As the result, permutations of lengths $2n+1$ and $2n$, respectively, are obtained; we use $\sigma$ to denote such a permutation. We claim that in either case, $\sigma$ is also $(k,\ell)$-crucial.  Indeed, non-extendability of $\sigma$ to the right follows from the fact that the black elements form, in the reduced form, the $(k,\ell)$-crucial permutation $\pi$ (instead of an arithmetic subsequence  with a difference $d$ in $\pi$ we deal with an arithmetic subsequence with the difference $2d$ in $\sigma$). On the other hand, $\sigma$ is $(k,\ell)$-anti-monotone since the longest monotone arithmetic progression with an odd difference $d$ is of length 2, while for  even $d$ any arithmetic subsequence is mono-chromatic (either entirely white or black) and we can use the properties of $\pi$ and $\tau$.     

\section{$(k,\ell)$-bicrucial permutations}\label{kl-bicrucial}

\begin{figure}[ht]
\begin{center}
\includegraphics[scale=0.5]{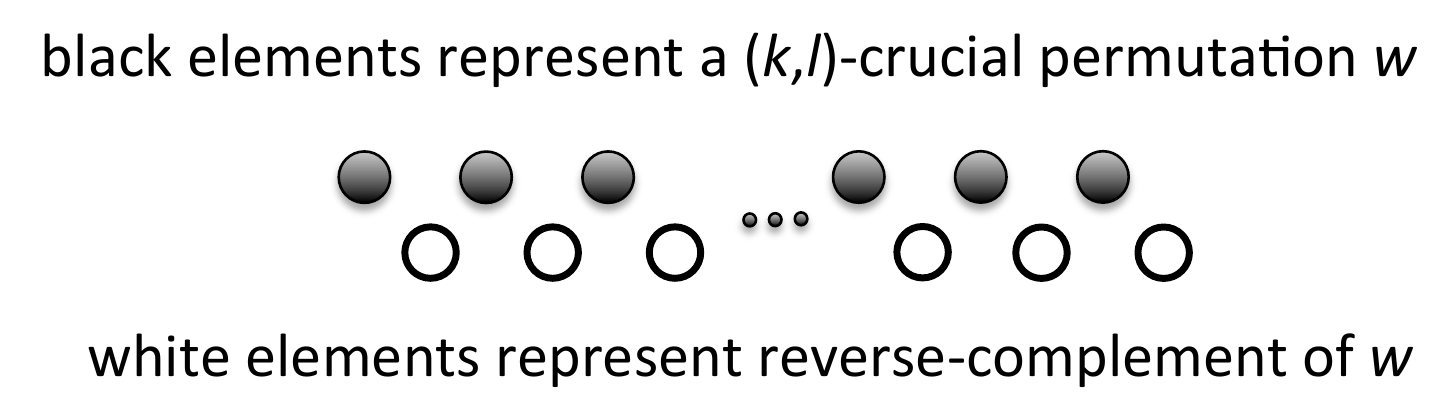}
\end{center}
\caption{A constructions of a $(k,\ell)$-bicrucial permutation of even length based on a $(k,\ell)$-crucial permutation.}\label{longBicrucial} 
\end{figure}

Based on the fact that $(k,\ell)$-crucial permutations exist for any $k,\ell\geq 3$, we can prove that $(k,\ell)$-bicrucial permutations exist for any $k,\ell\geq 3$. Indeed, let $w$ be a $(k,\ell)$-crucial permutation of length $n$, that is, $w$ is a $(k,\ell)$-anti-monotone permutation that is not expendable to the right. Then applying reverse complement to $w$, we will obtain a  $(k,\ell)$-anti-monotone permutation $w'$ that is not expendable to the left. Letting $w$ form the largest elements in a permutation of length $2n$ occupying every other position (as shown schematically in Figure \ref{longBicrucial} by black circles) and $w'$ form the remain elements (as shown schematically in Figure \ref{longBicrucial} by white circles), we obtain a $(k,\ell)$-bicrucial permutation.  Based on this construction, the minimal length of a $(k,\ell)$-bicrucial permutation is $\leq 2\max (k,\ell)(\min(k,\ell)-1)=2m(k,\ell)$. 

The length of a permutation  obtained in Figure \ref{longBicrucial} is even, but to get longer $(k,\ell)$-bicrucial permutations of odd lengths, we can take any $(k,\ell)$-bicrucial permutation of even length and shuffle it with a $(3,3)$-anti-monotone permutation as shown schematically to the left in Figure \ref{longCrucial} (one should ignore all the words in that picture). It is not difficult to check that this construction works.

\begin{figure}[ht]
\begin{center}
\includegraphics[scale=0.5]{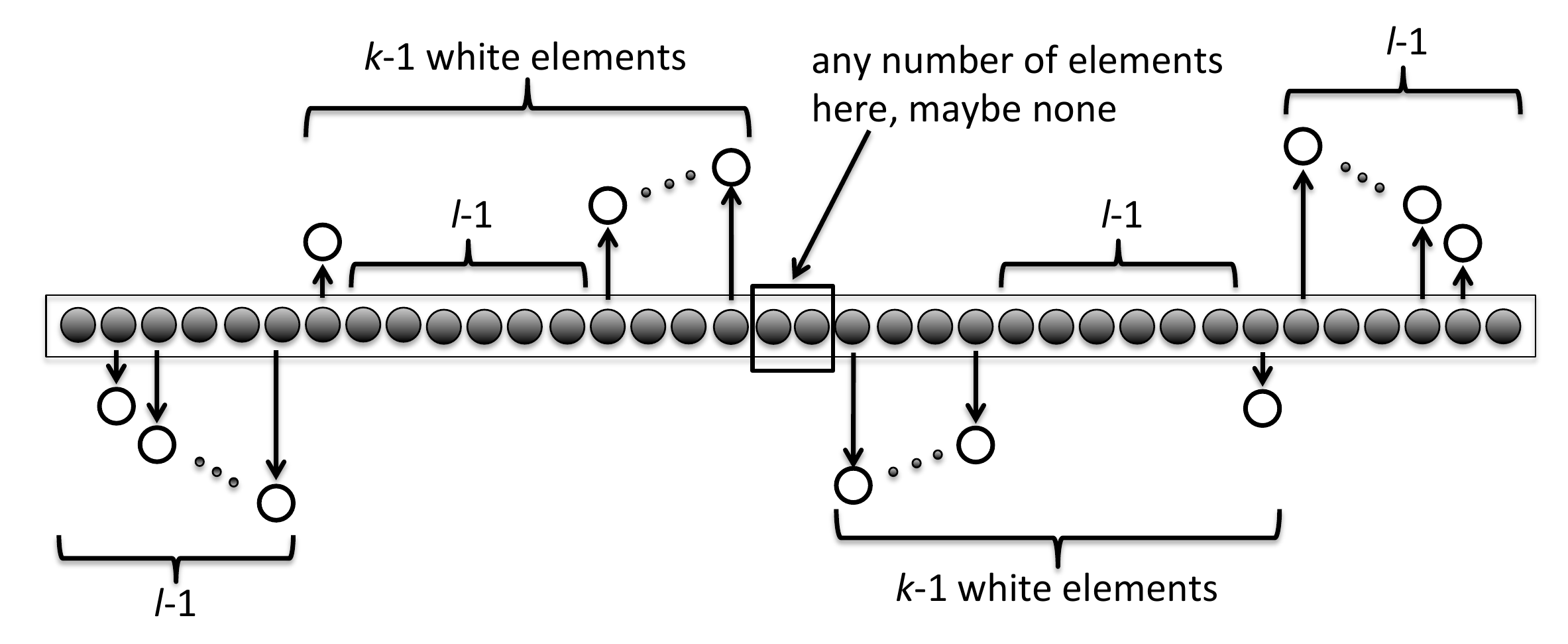}
\end{center}
\caption{A construction of $(k,\ell)$-bicrucial permutations, $k,\ell\geq 5$.}\label{constrBiCrucial} 
\end{figure}

An alternative construction of arbitrary long $(k,\ell)$-bicrucial permutations of length $2\ell(k-1)$ or more for $k,\ell\geq 5$, is given in Figure \ref{constrBiCrucial}; one can notice similarity with Figure \ref{mainConstruction} which gave us the idea of the construction. We note that there is some freedom in the construction in Figure \ref{constrBiCrucial}, namely the relative order of large elements (above the black elements) to the right and to the left of the middle is arbitrary; the same applies to the small elements (below the black elements) -- the relative order of those elements to the right and to the left of the middle is irrelevant.  A justification of the fact that the construction in Figure \ref{constrBiCrucial} works is similar to our justification of Figure \ref{mainConstruction} : non-extendability to the right and to the left is straightforward, while any arithmetic monotone subsequences of length 4 or more can have at most one large white element and one small white element leading to the fact that a monotone subsequence of the permutation in   Figure \ref{constrBiCrucial} are of length at most 4 (recall that the black elements form a (3,3)-anti-monotone permutation.

\begin{figure}[ht]
\begin{center}
\includegraphics[scale=0.5]{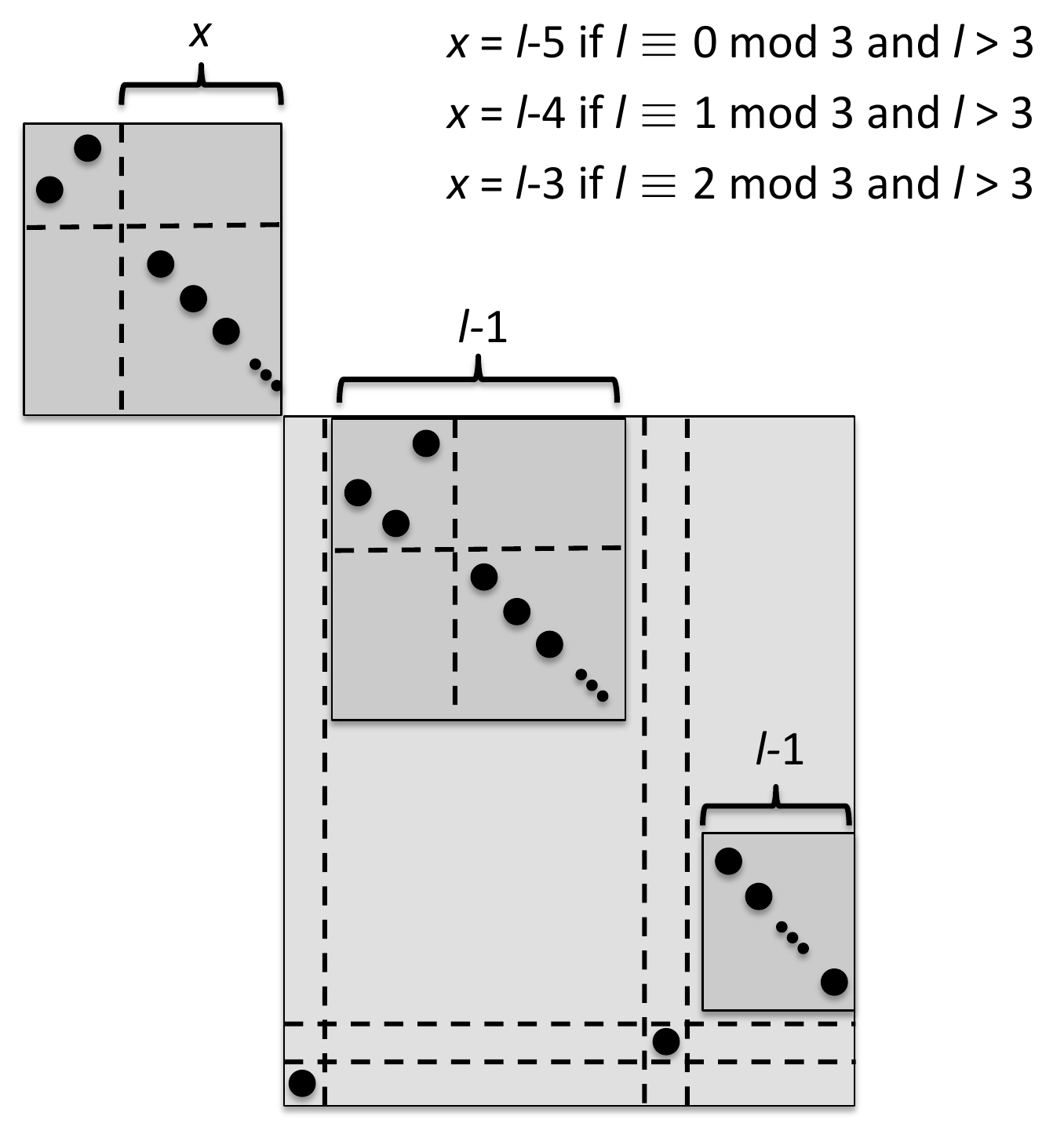}
\end{center}
\caption{A construction of $(3,\ell)$-bicrucial permutations, $\ell\geq 4$.}\label{3l-bicrucial} 
\end{figure}

Speaking of another alternative construction, which makes the length of a $(k,\ell)$-bicrucial permutation for the respective values of $k$ and $\ell$ smaller (and thus is a candidate for the minimal possible such length), in Figure \ref{3l-bicrucial} we provide a construction of $(3,\ell)$-bicrucial permutations where $\ell\geq 4$. The length of such permutations depends on which congruence class mod 3 $\ell$ belongs to. It is $3\ell-3+i$ if $\ell \equiv i\mod 3$, $i=0,1,2$. For example, for $\ell=4,5,6$ the construction gives the permutations $(10)(11)917682543$,  $(13)(14)(12)(11)198(10)726543$ and $(14)(15)(13) 1 (11)(10)(12) 98 2 76543$, respectively. 

To justify the construction in Figure \ref{3l-bicrucial} giving a permutation $\pi$, we notice that it is not extendable to the right (because it is based on the construction of a $(3,\ell)$-crucial permutation in Figure \ref{3l-construction}) and it is not extendable to the left (to avoid arithmetically 123 a prospective extension must be larger than the leftmost element of $\pi$ which will lead  to an arithmetic occurrence of the pattern $\ell(\ell-1)\ldots 1$ with the difference $d=3$; to see the last claim, one has to go through three cases depending on congruence class of $\ell$ modulo 3). It remains to justify that $\pi$ is $(3,\ell)$-anti-monotone. 

It is easy to see that $\pi$ does not contain an arithmetic occurrence of $123$. As for arithmetic occurrences of the pattern $\ell(\ell-1)\ldots 1$, clearly no occurrence exists for $d=1$, and no occurrence exists for $d=2$ which follows from the fact that $2\ell$ rightmost elements form a $(3,\ell)$-crucial permutation presented in Figure \ref{3l-construction}, while any other occurrence is not possible because of the element 1 and three elements immediately to the right of it. Also, no arithmetic occurrence of $\ell(\ell-1)\ldots 1$ is possible for the difference $d\geq 4$ because $\pi$ does not have enough elements. It only remains to show that no arithmetic occurrence of $\ell(\ell-1)\ldots 1$ is possible for the difference $d=3$ which we do separately in each of the following three cases:

\begin{itemize}
\item For $\ell\equiv 0\mod 3$, the only non-trivial case (giving the longest possible arithmetic occurrence of the pattern $m(m-1)\ldots1$ for some $m$) is starting with the rightmost element and jumping to the left with difference $d=3$, that is, taking from the right the elements in the 1st position, 4th, 7th, etc, $(3\ell-5)$th position. However, $m=\ell-1$ in this case, not $\ell$, and no occurrence of $\ell(\ell-1)\ldots 1$ is possible. 
\item For $\ell\equiv 1\mod 3$, the only non-trivial case (giving the longest possible arithmetic occurrence of the pattern $m(m-1)\ldots1$ for some $m$) is starting with the 3rd element from the right and jumping to the left with difference $d=3$, that is, taking from the right the elements in the 3rd position, 6th, 9th, etc, $(3\ell-4)$th position. However, $m=\ell-1$ in this case, not $\ell$, and no occurrence of $\ell(\ell-1)\ldots 1$ is possible. 
\item For $\ell\equiv 2\mod 3$, the only non-trivial case (giving the longest possible arithmetic occurrence of the pattern $m(m-1)\ldots1$ for some $m$) is starting with the 3rd element from the right and jumping to the left with difference $d=3$, that is, taking from the right the elements in the 3rd position, 6th, 9th, etc, $(3\ell-3)$th position. However, again $m=\ell-1$ in this case, not $\ell$, and no occurrence of $\ell(\ell-1)\ldots 1$ is possible. 
\end{itemize}

\section{Open problems}\label{final}

One can check by computer that no (3,3)-crucial permutation exists of length 9, and thus all (3,3)-crucial permutations of length 8 are (3,3)-bicrucial. This led us to the following question.   

\begin{problem} Classify lengths for which $(k,\ell)$-crucial and $(k,\ell)$-bicrucial permutations exist for $k,\ell\geq 3$. Partial classification follows from the results in this paper. \end{problem}

We have that the length of a minimal (3,3)-bicrucial permutation is 8 which is less than $2m(3,3)=12$. Also, for $\ell\geq 4$, Figure \ref{3l-bicrucial} gives examples of $(3,\ell)$-bicrucial permutations whose length is less than $2m(3,\ell)$. This led us to the following question.

\begin{problem} What is the minimal length of a $(k,\ell)$-bicrucial permutation?\end{problem}

Finally, a rich research direction is to study (bi)crucial permutations with respect to other sets of (arithmetic) prohibitions. We believe that this direction will bring many interesting and elegant results.

\end{document}